\documentclass[a4paper,11pt]{article}
\usepackage{amsmath,amsfonts,mathrsfs,amssymb,color}
\usepackage{indentfirst}
\usepackage{amsthm,mathtools}
\usepackage{microtype}
\usepackage{graphicx}
\usepackage{placeins}
\usepackage[hidelinks]{hyperref}
\hypersetup{
  pdfauthor={Ao Xiao},
  pdftitle={Exact multiplicity for a Minkowski-curvature Neumann problem with a step weight},
  pdfsubject={Positive solutions, exact multiplicity, and nondegenerate folds},
  pdfkeywords={Minkowski curvature, Neumann problem, indefinite weight, exact multiplicity, shooting method}
}

\numberwithin{equation}{section}

\newtheoremstyle{xiaopaper}
  {2mm}{2mm}{\normalfont}{0pt}{\bfseries}{.}{.75em}{}
\theoremstyle{xiaopaper}
\newtheorem{theorem}{Theorem}[section]
\newtheorem{proposition}[theorem]{Proposition}
\newtheorem{lemma}[theorem]{Lemma}
\newtheorem{corollary}[theorem]{Corollary}
\newtheorem{remark}[theorem]{Remark}

\makeatletter
\renewenvironment{proof}[1][\proofname]{%
  \par\pushQED{\qed}\normalfont
  \topsep2mm\relax
  \trivlist
  \item[\hskip\labelsep\bfseries #1\@addpunct{.}]\ignorespaces
}{%
  \popQED\endtrivlist\@endpefalse
}
\renewenvironment{thebibliography}[1]{%
  \par\addvspace{3mm}%
  {\centering\bfseries REFERENCES\par}\nopagebreak
  \@mkboth{\MakeUppercase\refname}{\MakeUppercase\refname}%
  \list{\@biblabel{\@arabic\c@enumiv}}{%
    \setlength{\leftmargin}{2.5em}%
    \setlength{\labelwidth}{2em}%
    \setlength{\labelsep}{0.5em}%
    \setlength{\itemindent}{0pt}%
    \setlength{\listparindent}{0pt}%
    \usecounter{enumiv}%
    \let\p@enumiv\@empty
    \renewcommand{\theenumiv}{\@arabic\c@enumiv}%
  }%
  \emergencystretch=3em
  \clubpenalty4000\widowpenalty4000\sfcode`.\@m
}{%
  \def\@noitemerr{\@latex@warning{Empty `thebibliography' environment}}%
  \endlist
}
\makeatother

\newcommand{\dd}{\,\mathrm d}
\newcommand{\asinh}{\operatorname{arsinh}}
\newcommand{\sgn}{\operatorname{sgn}}
\newcommand{\norm}[1]{\left\lVert#1\right\rVert}
\newcommand{\mcD}{\mathcal D}
\newcommand{\mcF}{\mathcal F}
\newcommand{\mcH}{\mathcal H}
\newcommand{\mcR}{\mathcal R}
\allowdisplaybreaks[2]

\title{Exact multiplicity for a Minkowski-curvature Neumann problem with a step weight}
\author{Ao Xiao$^{\ast}$\\
  {\small Department of Mathematics, Northwest Normal University,
   Lanzhou 730070, P R China}}
\date{}

\begin{document}

\maketitle
\footnote[0]{E-mail: \href{mailto:xiaoao2math@163.com}{xiaoao2math@163.com}.
  $^{\ast}$Corresponding author: Ao Xiao.}

\begin{abstract}
The DEG1 challenge asks whether positive solutions of a
super-exponential Minkowski-curvature equation are unique for step
weights with a single sign change and negative mean. In this paper,
we give a Neumann counterexample by studying
\[
\left(\frac{u'}{\sqrt{1-(u')^2}}\right)'
  +\lambda a_{\mu,h}(t)(e^{u^2}-1)=0,\qquad
u'(0)=u'(2h)=0,
\]
where $a_{\mu,h}=1$ on $(0,h)$ and $a_{\mu,h}=-\mu$ on $(h,2h)$,
with $\mu>1$. Using the shooting method, we prove that the positive
solution curve has exactly two nondegenerate folds for $(\lambda,h)$
in a nonempty open region. As $\mu$ increases through the corresponding
critical values, the number of positive solutions is successively
$1$, $2$, $3$, $2$, and $1$.
\end{abstract}


{\small\bfseries MR(2020)\quad 34B18, 34B15, 34C23.}

\baselineskip 20pt

\section{Introduction and main result}\label{sec:introduction}

The DEG1 challenge \cite{DEG1Challenge} asks whether the equation
\[
\left(\frac{u'}{\sqrt{1-(u')^2}}\right)'
  +a(t)\bigl(e^{u^p}-1\bigr)=0,\qquad p>1,
\]
has a unique positive solution under Neumann or periodic boundary
conditions when $a$ is a step function with one sign change and negative
mean. Boscaggin, Feltrin, and Zanolin \cite[Corollary~2.1]{BFZ2023}
proved existence for both boundary conditions under hypotheses that
cover these weights. They also noted that numerical experiments
suggested uniqueness in the super-exponential case
\cite[Section~4.2]{BFZ2021}.

We give a Neumann counterexample with $p=2$ and two steps of equal length.
Consider
\begin{equation}\label{eq:problem}
\begin{cases}
\displaystyle
\left(\frac{u'}{\sqrt{1-(u')^2}}\right)'
  +\lambda a_{\mu,h}(t)\bigl(e^{u^2}-1\bigr)=0,
  &0<t<2h,\\[6pt]
u'(0)=u'(2h)=0,
\end{cases}
\end{equation}
where $\lambda,h>0$, $\mu>1$, and
\begin{equation}\label{eq:weight}
a_{\mu,h}(t)=
\begin{cases}
1,&0<t<h,\\
-\mu,&h<t<2h.
\end{cases}
\end{equation}
A positive solution is a function $u\in C^1([0,2h])$ with $u>0$ and
$|u'|<1$ on $[0,2h]$, whose flux $u'/\sqrt{1-(u')^2}$ is absolutely
continuous and which satisfies \eqref{eq:problem} almost everywhere.
We denote the number of distinct positive solutions by $N(\lambda,h,\mu)$.

Multiplicity results for indefinite Minkowski-curvature equations
include pairs of positive radial Neumann solutions for power
nonlinearities \cite{BF2020} and positive periodic and subharmonic
solutions \cite{BFPeriodic2020}. Shooting methods have also been
used to study multiplicity for radial Neumann problems \cite{BCN2020}.
Ma, Su, and Zhao \cite{MSZ2023} studied the Neumann problem
\[
-\left(\frac{u'}{\sqrt{1-(u')^2}}\right)'=\lambda a(t)f(u),
\quad 0<t<T,\qquad u'(0)=u'(T)=0,
\]
with an integrable weight $a$ allowed to change sign. Under suitable
growth conditions on $f$, they obtained at least two positive solutions
for sufficiently large $\lambda$. For the same form of equation
on $(0,1)$, Su, Ma, and Yang \cite{SMY2025} fixed a continuous
sign-changing weight $a$ and used bifurcation methods to establish
the existence of $S$-shaped connected components of positive solutions
as the overall multiplier $\lambda$ varies. In \eqref{eq:problem} we
fix the positive coefficient $\lambda$ and vary the negative one
through the ratio $\mu$. This lets us study how the number of positive
solutions changes with the balance between the two steps.

In the semilinear setting, Feltrin and Sovrano \cite[Theorem~1.1]{FS2018}
used shooting to obtain at least three positive Neumann solutions for
$-u''=(\lambda a^+-\mu a^-)g(u)$ under their assumptions. Feltrin,
Sovrano, and Tellini \cite[Section~1]{FST2022} studied solution counts
for $-u''=a_{\lambda,\mu}(t)u^2(1-u)$ with Neumann conditions and a
step weight having two positive intervals separated by a negative one.
Exact multiplicity and bifurcation curves for autonomous
Minkowski-curvature Dirichlet problems were studied in
\cite{Huang2018,Huang2019}. In particular, Huang and Wang
\cite[Theorem~4]{HW2025} studied the autonomous perturbed Gelfand problem
\[
-\left(\frac{u'}{\sqrt{1-(u')^2}}\right)'
=\lambda\exp\!\left(\frac{au}{a+u}\right),
\quad -d<t<d,\qquad u(-d)=u(d)=0.
\]
They used a time map relating the maximum of a positive solution
to the half-length $d$ to classify the shape of the
$(\lambda,\|u\|_\infty)$ curve as $a$ and $d$ vary. In particular,
the curve has exactly two turning points when $a$ is above a
threshold and $d$ is sufficiently large. For nonautonomous Dirichlet problems,
$S$-shaped connected components and more complicated bifurcation
curves have been obtained in \cite{MX2018,LSY2024}.

For $(\lambda,h)$ in the open region of Theorem~\ref{thm:main},
matching the trajectories across the jump in the weight gives exactly
one $\mu>1$ for each initial value $s=u(0)>0$. The resulting graph
parametrizes all positive solutions of \eqref{eq:problem}. A uniform
$C^2$ expansion on bounded amplitude intervals gives two nondegenerate
folds. A variational
estimate proves strict monotonicity at large amplitudes, excluding
further folds on the unbounded part of the graph. Counting intersections
with the lines of constant $\mu$ then gives the exact number of positive
solutions for every $\mu>1$ in this regime.

\begin{theorem}\label{thm:main}
There exist constants $L_*>0$ and $h_*>0$ such that, for every
\begin{equation}\label{eq:region}
(\lambda,h)\in\mcD
:=\bigl\{(\lambda,h)\in(0,\infty)^2:
       0<h<h_*,\quad L_*<\lambda h<2L_*\bigr\},
\end{equation}
there are unique critical values $1<\mu_-(\lambda,h)<\mu_+(\lambda,h)<\infty$
such that, for every $\mu>1$,
\begin{equation}\label{eq:count}
N(\lambda,h,\mu)=
\begin{cases}
1,&1<\mu<\mu_-(\lambda,h),\\
2,&\mu=\mu_-(\lambda,h),\\
3,&\mu_-(\lambda,h)<\mu<\mu_+(\lambda,h),\\
2,&\mu=\mu_+(\lambda,h),\\
1,&\mu>\mu_+(\lambda,h).
\end{cases}
\end{equation}
The functions $\mu_-$ and $\mu_+$ are smooth on $\mcD$.
For each fixed $(\lambda,h)\in\mcD$, the set of positive solutions is
a single smooth curve parametrized by $s=u(0)>0$. Its projection onto
the $\mu$-axis has exactly two critical points, both nondegenerate.
At each critical value precisely one solution is degenerate, with a
one-dimensional kernel for the linearized Neumann problem. Every other
positive solution is nondegenerate.
\end{theorem}

For any $(\lambda,h)\in\mcD$ and
$\mu\in(\mu_-(\lambda,h),\mu_+(\lambda,h))$, the weight
$\lambda a_{\mu,h}$ has one sign change and integral
$\lambda h(1-\mu)<0$. The three positive solutions therefore disprove
uniqueness in the Neumann part of the DEG1 challenge. Uniqueness is
recovered for $\mu>\mu_+(\lambda,h)$.

Section~\ref{sec:shooting} constructs the global shooting curve.
Section~\ref{sec:proof} analyzes its shape and proves the main theorem.
Section~\ref{sec:asymptotics} derives the critical-value asymptotics
and illustrates the curve numerically.

\section{The global shooting curve}\label{sec:shooting}

Throughout the paper, set
\begin{equation}\label{eq:definitions}
\begin{gathered}
g(u)=e^{u^2}-1,\qquad G(u)=\int_0^u g(r)\dd r,\\
\phi(p)=\frac{p}{\sqrt{1-p^2}},\qquad
\psi(v)=\frac{v}{\sqrt{1+v^2}},\qquad
K(v)=\sqrt{1+v^2}-1.
\end{gathered}
\end{equation}
Thus $\psi=\phi^{-1}$ and $K'=\psi$. With $v=-\phi(u')$,
problem \eqref{eq:problem} becomes
\begin{equation}\label{eq:system}
u'=-\psi(v),\qquad
v'=\begin{cases}
\lambda g(u),&0<t<h,\\
-\lambda\mu g(u),&h<t<2h,
\end{cases}
\qquad v(0)=v(2h)=0.
\end{equation}
For each finite initial condition, $|u'|<1$ bounds $u$ on $[0,2h]$,
and the equation for $v$ then bounds $v'$. Thus neither state variable
can escape to infinity in finite time. Applying the continuation
criterion \cite[Corollary~2.16]{Teschl2012} to each autonomous system
gives existence throughout $[0,2h]$ by concatenation at $h$. Both
vector fields are smooth, so their flows depend smoothly on the
initial data and parameters by \cite[Theorem~2.11]{Teschl2012}.
The composition preserves this smooth dependence.

\begin{lemma}\label{lem:decreasing}
Every positive solution of \eqref{eq:problem} satisfies
\begin{equation}\label{eq:positive-momentum}
v(t)>0\quad(0<t<2h),\qquad u'(t)<0\quad(0<t<2h).
\end{equation}
Moreover,
\begin{equation}\label{eq:balance}
\mu=
\frac{\displaystyle\int_0^h g(u(t))\dd t}
     {\displaystyle\int_h^{2h}g(u(t))\dd t}>1.
\end{equation}
\end{lemma}

\begin{proof}
The sign assertions in \eqref{eq:positive-momentum} follow from
\cite[Section~2.1, Eqs.~(2.5)--(2.6)]{BFZ2021} with
$(x,y)=(u,-v)$, $(a_+,a_-)=(\lambda,\lambda\mu)$, and $(\tau,T)=(h,2h)$.
Integrating the equation for $v$ gives the ratio in \eqref{eq:balance}.
It is greater than one because $g$ is strictly increasing and
$u(t)>u(t+h)$ for $0\le t\le h$.
\end{proof}

We now fix $L=\lambda h$. The following condition ensures that every positive initial
value can be used for shooting:
\begin{equation}\label{eq:left-condition}
Lh^2e^{h^2}<1.
\end{equation}

\begin{lemma}\label{lem:left-positive}
Assume \eqref{eq:left-condition}. For every $s>0$, the solution of
the first-half system in \eqref{eq:system} with $u(0)=s$ and $v(0)=0$
is positive on $[0,h]$. In particular,
\begin{equation}\label{eq:cq}
c(s):=u(h)>0,\qquad q(s):=v(h)>0.
\end{equation}
\end{lemma}

\begin{proof}
If $s\ge h$, then $u(t)>s-t\ge0$ for $0<t\le h$, since $|u'|<1$.
Suppose that $0<s<h$ and that $u$ has a first zero in $[0,h]$.
Before this zero, $0<u\le s$ and $0\le v(t)\le\lambda h g(s)=Lg(s)$.
Consequently,
\[
s-u(t)\le hLg(s)
\le hLs^2e^{h^2}
<\bigl(Lh^2e^{h^2}\bigr)s<s,
\]
which is impossible at the first zero.
\end{proof}

\begin{proposition}\label{prop:global-curve}
Assume \eqref{eq:left-condition}. For each $s>0$, there is exactly one
$\mu=\mu_{L,h}(s)>0$ for which the solution of \eqref{eq:system} with
$u(0)=s$ is positive and satisfies $v(2h)=0$. The function
$\mu_{L,h}:(0,\infty)\to(1,\infty)$ is smooth, and its graph
parametrizes all positive solutions.

For fixed $L,h$, let
\begin{equation}\label{eq:shooting-F}
F(s,\mu)=v(2h;s,\mu),
\end{equation}
where the dependence on $L,h$ is suppressed. At every positive solution,
\begin{equation}\label{eq:Fmu}
F_\mu(s,\mu)<0.
\end{equation}
\end{proposition}

\begin{proof}
Fix $s>0$ and use \eqref{eq:cq} as initial data for the second
half of the trajectory. In this proof only, time on that half is
measured from zero. We use the energy identity and time quadrature
for the negative step from \cite[Section~2.1, Eqs.~(2.7)--(2.8)]{BFZ2021}.
In our notation, the conserved quantity is
\begin{equation}\label{eq:energy}
K(v)-\lambda\mu G(u)=K(q)-\lambda\mu G(c).
\end{equation}
A positive stopping point $(b,0)$ therefore satisfies
\begin{equation}\label{eq:b}
G(b)=G(c)-\frac{K(q)}{\lambda\mu}.
\end{equation}
Such a point exists precisely when
\begin{equation}\label{eq:mu0}
\mu>\mu_0(s):=\frac{K(q)}{\lambda G(c)}.
\end{equation}
In that case it is unique and belongs to $(0,c)$.
The trajectory has $u'<0$ and $v'<0$ before reaching $(b,0)$.
Solving \eqref{eq:energy} for $v$ gives its stopping time:
\begin{equation}\label{eq:time-map}
T_s(\mu)=\int_b^c
\frac{1+\lambda\mu\bigl(G(u)-G(b)\bigr)}
{\sqrt{\lambda\mu\bigl(G(u)-G(b)\bigr)
\bigl[2+\lambda\mu\bigl(G(u)-G(b)\bigr)\bigr]}}\dd u.
\end{equation}
The singularity at $u=b>0$ is of order $(u-b)^{-1/2}$, so this time
is finite. The stopping event is transverse because
$v'(T_s)=-\lambda\mu g(b)<0$.

To prove strict monotonicity in $\mu$, write $z=-v$ on the second
half-interval. Then $u'=\psi(z)$ and $z'=\lambda\mu g(u)$, with
$(u(0),z(0))=(c,-q)$.
The derivatives $X=\partial_\mu u$ and $Y=\partial_\mu z$ satisfy
\begin{equation}\label{eq:mu-variation}
X'=\psi'(z)Y,\qquad
Y'=\lambda g(u)+\lambda\mu g'(u)X,
\qquad X(0)=Y(0)=0.
\end{equation}
On the positive trajectory, this is a cooperative linear system with
a strictly positive forcing in its second equation. Hence $X(t)>0$
and $Y(t)>0$ for $0<t\le T_s(\mu)$. Differentiating
$z(T_s(\mu),\mu)=0$ yields
\begin{equation}\label{eq:Tmu}
T_s'(\mu)=-\frac{Y(T_s(\mu))}{\lambda\mu g(b)}<0.
\end{equation}

As $\mu\downarrow\mu_0(s)$, equation \eqref{eq:b} gives $b\downarrow0$.
The integrand in \eqref{eq:time-map}, extended by zero to $(0,b)$,
converges pointwise on $(0,c)$ to
\[
\frac{1+\lambda\mu_0G(u)}
{\sqrt{\lambda\mu_0G(u)\bigl(2+\lambda\mu_0G(u)\bigr)}}.
\]
Since $G(u)\sim u^3/3$ at zero, this limit is asymptotic to a positive
constant times $u^{-3/2}$. Fatou's lemma gives
$T_s(\mu)\to\infty$.
As $\mu\to\infty$, we have $b\to c$, and integration of
$-v'=\lambda\mu g(u)\ge\lambda\mu g(b)$ gives
$T_s(\mu)\le q/(\lambda\mu g(b))\to0$.
Thus there is exactly one $\mu>\mu_0(s)$ with $T_s(\mu)=h$.
The resulting solution is positive on the whole interval, and
Lemma~\ref{lem:decreasing} shows that $\mu>1$.
Conversely, every positive Neumann solution must stop in this way.

At a matched trajectory, \eqref{eq:mu-variation}, evaluated at the
fixed right-hand time $h$, gives $F_\mu=-Y(h)<0$.
The ordinary implicit function theorem therefore gives smooth dependence
of $\mu_{L,h}(s)$ on $s$, and also on $L,h$ as long as
\eqref{eq:left-condition} holds. Uniqueness identifies the local
implicit functions with the global curve just constructed.
\end{proof}

\begin{lemma}\label{lem:infinity}
Under \eqref{eq:left-condition},
\begin{equation}\label{eq:infinity}
\lim_{s\to\infty}\mu_{L,h}(s)=\infty.
\end{equation}
\end{lemma}

\begin{proof}
For $s>h$, the first half of the trajectory satisfies
$u(t)>s-h$ and $v(t)\ge A_st$, where $A_s:=\lambda g(s-h)$.
It follows that
\begin{align}
c(s)
&\le s-\int_0^h\psi(A_st)\dd t
 =s-\frac{\sqrt{1+A_s^2h^2}-1}{A_s}\notag\\
&\le s-h+\frac1{A_s}.
\label{eq:c-large}
\end{align}
For all sufficiently large $s$, this gives $c(s)\le s-3h/4$.
On $[0,h/2]$ we have $u(t)\ge s-h/2$, whereas on $[h,2h]$
we have $u(t)\le c(s)$. Applying \eqref{eq:balance}, we obtain
\[
\mu_{L,h}(s)
\ge\frac{g(s-h/2)}{2g(s-3h/4)}\longrightarrow\infty.
\]
The last limit follows from
$(s-h/2)^2-(s-3h/4)^2=hs/2-5h^2/16$.
\end{proof}

\section{Proof of the main theorem}\label{sec:proof}

By Proposition~\ref{prop:global-curve}, for fixed $L=\lambda h$ and $h$,
counting positive solutions amounts to counting the initial values
$s>0$ for which $\mu_{L,h}(s)=\mu$. We first analyze the leading term of
the shooting curve on bounded intervals. We then prove that the two
critical points persist and that the curve is strictly increasing
for $s\ge2$, so these are its only turning points. The two endpoint
limits and the three monotone branches give the exact count, and
the linearized shooting equation identifies the degenerate solutions.

\subsection{The reduced shooting curve}

\subsubsection{A uniform expansion on bounded amplitude intervals}\label{sec:expansion}

Define, for $s>0$ and $k>0$,
\begin{equation}\label{eq:ell-W}
\ell(s)=\frac{g'(s)}{g(s)}=\frac{2s}{1-e^{-s^2}},\qquad
W(k)=\frac{k\sqrt{1+k^2}-\asinh k}{k^2},
\end{equation}
and set $W(0)=0$. An integral representation that will be useful below is
\begin{equation}\label{eq:W-integral}
W(k)=2\int_0^1r\psi(kr)\dd r
    =2k\int_0^1\frac{r^2}{\sqrt{1+k^2r^2}}\dd r.
\end{equation}
The reduced shooting function is
\begin{equation}\label{eq:M}
M_L(s)=\ell(s)W\bigl(Lg(s)\bigr),\qquad s>0.
\end{equation}

\begin{proposition}\label{prop:expansion}
Let $S>0$ and let $J\Subset(0,\infty)$ be a compact interval.
For sufficiently small $h>0$, uniformly for $L\in J$, the shooting
function extends smoothly to $s=0$, with $\mu_{L,h}(0)=1$.
There is a constant $C_{J,S}$ such that
\begin{equation}\label{eq:C2-expansion}
\sup_{L\in J}
\norm{\mu_{L,h}-1-hM_L}_{C^2([0,S])}
\le C_{J,S}h^2.
\end{equation}
Here $M_L$ has a smooth extension to $s=0$, and
\begin{equation}\label{eq:M-zero}
M_L(0)=0,\qquad M_L'(0)=\frac{4L}{3}.
\end{equation}
The function $(\mu_{L,h}-1)/h$ itself extends smoothly to $h=0$,
where its value is $M_L$.
\end{proposition}

\begin{proof}
Put $x=t/h$ and write $u(hx)=sy(x)$ and $v(hx)=Lg(s)z(x)$.
For $s>0$, the normalized equations are
\begin{equation}\label{eq:normalized}
\begin{aligned}
y_x&=-h\frac{Lg(s)}s\,
       \frac{z}{\sqrt{1+[Lg(s)z]^2}},\\
z_x&=\widehat a_\mu(x)\frac{g(sy)}{g(s)},
\qquad y(0)=1,\quad z(0)=0,
\end{aligned}
\end{equation}
where $\widehat a_\mu=1$ on $(0,1)$ and
$\widehat a_\mu=-\mu$ on $(1,2)$.
Let $E(r)=(e^r-1)/r$ for $r\ne0$, with the smooth positive extension
$E(0)=1$.
The identities
\begin{equation}\label{eq:removable}
\frac{g(s)}s=sE(s^2),\qquad
\frac{g(sy)}{g(s)}=y^2\frac{E(s^2y^2)}{E(s^2)}
\end{equation}
remove the apparent singularities at $s=0$.
Thus the two vector fields in \eqref{eq:normalized} are smooth
in $(s,L,h,\mu,y,z)$ near the compact parameter set under consideration.

At $h=0$ the solution is
\begin{equation}\label{eq:normalized-base}
y\equiv1,\qquad
z(x)=\begin{cases}x,&0\le x\le1,\\
1-\mu(x-1),&1\le x\le2.
\end{cases}
\end{equation}
Smooth dependence \cite[Theorem~2.11]{Teschl2012}, applied successively
on the two half-intervals, provides a common neighborhood of
$[0,S]\times J\times\{0\}\times\{1\}$ in which the flow exists up
to $x=2$ and depends smoothly on all parameters.
For its terminal value $\mcF(s,L,h,\mu):=z(2;s,L,h,\mu)$, we have
\begin{equation}\label{eq:normalized-IFT}
\mcF(s,L,0,\mu)=1-\mu,\qquad
\partial_\mu\mcF(s,L,0,1)=-1.
\end{equation}
The implicit function theorem and compactness therefore yield a single
smooth function $\mu=\Gamma(s,L,h)$, defined for
$s\in[0,S]$, $L\in J$, and $|h|$ small, with
$\Gamma(s,L,0)=1$.
At $s=0$, equations \eqref{eq:normalized}--\eqref{eq:removable} give
$y\equiv1$ for every such $h$, and hence $\Gamma(0,L,h)=1$.
For $s>0$ and small positive $h$, the normalized solution has $y>0$.
After also imposing \eqref{eq:left-condition},
Proposition~\ref{prop:global-curve} identifies
$\Gamma(s,L,h)$ with $\mu_{L,h}(s)$.

We compute the derivative with respect to $h$ at zero.
Fix $s>0$ and put $k=Lg(s)$. Along the implicit branch, let
$U_h(x)=u(hx)$ and $V_h(x)=v(hx)$, where the subscript denotes the
parameter, not differentiation. At $h=0$,
\[
U_0(x)=s,\qquad
V_0(x)=\begin{cases}kx,&0\le x\le1,\\k(2-x),&1\le x\le2.\end{cases}
\]
Since $(U_h)_x=-h\psi(V_h)$, its parameter derivative satisfies
\[
\dot U(x):=\left.\partial_hU_h(x)\right|_{h=0}
=-\int_0^x\psi(V_0(r))\dd r.
\]
Rescaling \eqref{eq:balance} to $[0,2]$ and differentiating gives
\begin{align*}
\partial_h\Gamma(s,L,0)
&=\ell(s)\left(\int_0^1\dot U(x)\dd x
                  -\int_1^2\dot U(x)\dd x\right)\\
&=2\ell(s)\int_0^1r\psi(kr)\dd r
=M_L(s).
\end{align*}
The second equality follows by changing the order of integration on
each half-interval.

The quotient $W(k)/k$ extends smoothly to zero and
$W(k)=2k/3+O(k^3)$. Writing $M_L(s)=Lg'(s)W(Lg(s))/(Lg(s))$
proves smooth extendibility at $s=0$ and \eqref{eq:M-zero}.
Taylor's formula now gives
\begin{equation}\label{eq:remainder}
\Gamma(s,L,h)=1+hM_L(s)+h^2\mcR(s,L,h),
\end{equation}
where
\[
\mcR(s,L,h)=\int_0^1(1-\theta)
          \partial_h^2\Gamma(s,L,\theta h)\dd\theta
\]
is smooth on a neighborhood of the compact set.
Its first two $s$-derivatives are uniformly bounded, proving
\eqref{eq:C2-expansion}. Equation \eqref{eq:remainder} also proves
the final assertion.
\end{proof}

\subsubsection{Two critical points of the reduced curve}\label{sec:reduced}

We first record the one-variable identities underlying the shape of
$M_L$. Set
\begin{equation}\label{eq:QR}
Q(k)=\frac{kW'(k)}{W(k)}\quad(k>0),\qquad
R(x)=\frac{(1+2x)e^{-x}-1}{2x}\quad(x>0),
\end{equation}
with $Q(0)=1$ and $R(0)=1/2$.

\begin{lemma}\label{lem:Q}
The function $Q$ is smooth at zero and strictly decreasing on
$(0,\infty)$. It satisfies
\begin{equation}\label{eq:Q-bounds}
\frac1{1+k^2}<Q(k)<1,\qquad
\lim_{k\to0}Q(k)=1,\qquad
\lim_{k\to\infty}Q(k)=0,
\end{equation}
and
\begin{equation}\label{eq:Q-identity}
kQ'(k)=(Q(k)+2)\left(\frac1{1+k^2}-Q(k)\right).
\end{equation}
In particular, there is a unique $\kappa>0$ such that
$Q(\kappa)=1/2$, and $Q'(\kappa)<0$.
\end{lemma}

\begin{proof}
Differentiating \eqref{eq:W-integral} yields
\[
Q(k)=
\frac{\displaystyle\int_0^1r^2(1+k^2r^2)^{-3/2}\dd r}
     {\displaystyle\int_0^1r^2(1+k^2r^2)^{-1/2}\dd r}.
\]
This is a weighted average of $(1+k^2r^2)^{-1}$, which proves the
strict inequalities in \eqref{eq:Q-bounds}.
The formula for $W$ also gives
\begin{equation}\label{eq:W-identity}
kW'(k)+2W(k)=\frac{2k}{\sqrt{1+k^2}}.
\end{equation}
Taking a logarithmic derivative of
$Q(k)+2=2k/(\sqrt{1+k^2}W(k))$ proves
\eqref{eq:Q-identity}; its right-hand side is negative.
The expansion of $W$ at zero proves smoothness and $Q(0)=1$.
Finally, dominated convergence in \eqref{eq:W-integral} gives
$W(k)\to1$ as $k\to\infty$. Equation \eqref{eq:W-identity} then
implies $Q(k)\to0$.
\end{proof}

The equation
\begin{equation}\label{eq:x0}
e^{x_0}=1+2x_0
\end{equation}
has exactly one positive root, and $1<x_0<2$. Indeed,
$e^x-1-2x$ is strictly convex, vanishes at zero, has negative
derivative there, and changes sign between $1$ and $2$.
Consequently,
\begin{equation}\label{eq:R-sign}
R(x)>0\quad(0\le x<x_0),\qquad
R(x)\le0\quad(x\ge x_0),\qquad R'(x_0)<0.
\end{equation}
A direct differentiation of $\ell$ gives
\begin{equation}\label{eq:ell-identity}
-\frac{\ell'(s)}{\ell(s)^2}=R(s^2).
\end{equation}
Combining this identity with \eqref{eq:M} yields
\begin{equation}\label{eq:M-derivative}
M_L'(s)=\ell(s)^2W\bigl(Lg(s)\bigr)
\left[Q\bigl(L(e^{s^2}-1)\bigr)-R(s^2)\right].
\end{equation}

\begin{proposition}\label{prop:two-critical}
There exists $L_0>0$ such that, for every $L>L_0$, the function
$M_L$ has exactly two critical points
\[
0<\alpha_L<\beta_L<\sqrt{x_0}<2.
\]
They are nondegenerate, with
\begin{equation}\label{eq:M-sign}
M_L'(s)\begin{cases}
>0,&0<s<\alpha_L,\\
<0,&\alpha_L<s<\beta_L,\\
>0,&s>\beta_L,
\end{cases}
\qquad
M_L''(\alpha_L)<0<M_L''(\beta_L).
\end{equation}
The functions $L\mapsto\alpha_L$ and $L\mapsto\beta_L$ are smooth
on $(L_0,\infty)$.
\end{proposition}

\begin{proof}
By \eqref{eq:M-derivative}, the critical points correspond to the
positive zeros of
\begin{equation}\label{eq:D-L}
D_L(x)=Q\bigl(L(e^x-1)\bigr)-R(x),\qquad x=s^2.
\end{equation}
There are no zeros for $x\ge x_0$, by \eqref{eq:R-sign}.

Choose $\varepsilon>0$ so small that $x_0-\varepsilon>0$ and
\begin{equation}\label{eq:R-negative-derivative}
R'(x)\le-c_1<0
\quad\text{for }x\in[x_0-\varepsilon,x_0].
\end{equation}
Next choose $0<\delta<x_0-\varepsilon$, and set
$c_0=\min_{[0,\delta]}R>0$.
By Lemma~\ref{lem:Q}, we can choose a fixed $K>\kappa$ with
$Q(K)<c_0$.

For the first zero, use the coordinate $k=L(e^x-1)$ and consider
\begin{equation}\label{eq:H-L}
H_L(k)=Q(k)-R\bigl(\log(1+k/L)\bigr),\qquad 0\le k\le K.
\end{equation}
As $L\to\infty$,
\[
H_L\longrightarrow Q-\frac12
\quad\text{in }C^1([0,K]).
\]
The limit has exactly one zero, $\kappa$, and its derivative there
is negative. To see that the same exact count holds for $H_L$, choose
a closed interval around $\kappa$ on which $Q'$ is bounded above
by a negative constant. Outside this interval,
$|Q-1/2|$ has a positive lower bound. The $C^1$ convergence gives
one strictly decreasing zero inside the interval and no zeros outside.
For sufficiently large $L$, the whole range $0\le k\le K$
corresponds to $0\le x<\delta$.
If $x\le\delta$ but $k=L(e^x-1)\ge K$, then
\begin{equation}\label{eq:small-gap}
D_L(x)\le Q(K)-c_0<0.
\end{equation}
Thus there is exactly one zero in $[0,\delta]$, and its sign change
is from positive to negative.

On $[\delta,x_0-\varepsilon]$, the function $R$ has a positive
minimum, while $Q(L(e^x-1))\to0$ uniformly. Hence $D_L<0$
throughout this interval for large $L$.

On $[x_0-\varepsilon,x_0]$, differentiation and
\eqref{eq:Q-identity} give, with $k=L(e^x-1)$,
\begin{equation}\label{eq:Q-x-derivative}
\frac{\mathrm d}{\mathrm dx}Q\bigl(L(e^x-1)\bigr)
=\frac{e^x}{e^x-1}(Q(k)+2)
 \left(\frac1{1+k^2}-Q(k)\right).
\end{equation}
The right-hand side tends uniformly to zero. Together with
\eqref{eq:R-negative-derivative}, this implies
$D_L'(x)\ge c_1/2>0$ for all sufficiently large $L$.
We already know that $D_L(x_0-\varepsilon)<0$, whereas
$D_L(x_0)=Q(L(e^{x_0}-1))>0$.
There is therefore exactly one zero in this last interval,
with sign change from negative to positive.

These intervals cover $[0,x_0]$, and both zeros are simple.
Since the prefactor in \eqref{eq:M-derivative} is positive for
$s>0$, the conclusions in \eqref{eq:M-sign} follow.
Smoothness of the two critical points is a consequence of the
implicit function theorem and their uniqueness.
\end{proof}

\subsection{The two folds of the shooting curve}

\subsubsection{Strict monotonicity at large amplitudes}\label{sec:tail}

The preceding expansion controls only bounded intervals of $s$.
We next establish a uniform estimate on the remaining half-line.

\begin{proposition}\label{prop:tail}
Let $J=[L_-,L_+]\Subset(0,\infty)$. There exists $h_J>0$ such that,
for $L\in J$ and $0<h<h_J$, the global shooting curve is defined and
\begin{equation}\label{eq:tail-positive}
\mu_{L,h}'(s)>0\qquad\text{for every }s\ge2.
\end{equation}
\end{proposition}

\begin{proof}
Choose $h_J\le1/4$ small enough that
\eqref{eq:left-condition} holds for every $L\in J$.
Fix $s\ge2$ and write $u(t;s),v(t;s)$ for the positive solution
on the shooting curve. Since $|u'|<1$,
\begin{equation}\label{eq:large-range}
u(t;s)>s-2h\ge\frac32\qquad(0\le t\le2h).
\end{equation}
The logarithmic derivative satisfies
\begin{equation}\label{eq:ell-positive-derivative}
\ell'(u)=\frac{2\bigl[1-(1+2u^2)e^{-u^2}\bigr]}
                  {(1-e^{-u^2})^2}>0
\qquad\text{for }u\ge\frac32,
\end{equation}
because $u^2>x_0$ in this range.

We differentiate along the shooting curve, keeping $L,h$ fixed:
\begin{equation}\label{eq:AB}
A(t)=\frac{\mathrm d}{\mathrm ds}u(t;s),\qquad
B(t)=\frac{\mathrm d}{\mathrm ds}v(t;s).
\end{equation}
On the first half-interval these are also the fixed-$\mu$ derivatives,
since the left equation does not involve $\mu$. They satisfy
\begin{equation}\label{eq:left-variation}
A'=-\psi'(v)B,\qquad B'=\lambda g'(u)A,
\qquad A(0)=1,\quad B(0)=0.
\end{equation}
Before a first possible zero of $A$, we have $0<A\le1$ and $B>0$.
Moreover, $u\le s$ and \eqref{eq:ell-positive-derivative} give
$B'=v'\ell(u)A\le\ell(s)v'$, hence $B(t)\le\ell(s)v(t)$.
On this half-interval $v'\ge\lambda g(s-h)$, so
\begin{align}
1-A(t)
&=\int_0^t\psi'(v(\tau))B(\tau)\dd\tau\notag\\
&\le\frac{\ell(s)}{\lambda g(s-h)}
   \int_0^{v(t)}\frac{r}{(1+r^2)^{3/2}}\dd r
<\frac{\ell(s)}{\lambda g(s-h)}.
\label{eq:A-bound}
\end{align}
The constant
\begin{equation}\label{eq:tail-constant}
C_\infty:=\sup_{s\ge2}\frac{\ell(s)}{g(s-1/4)}
\end{equation}
is finite: the ratio is continuous and tends to zero at infinity.
Since $\lambda=L/h$ and $h\le1/4$, we have
$\ell(s)/(\lambda g(s-h))\le hC_\infty/L_-$. Decrease $h_J$ so that
$hC_\infty/L_-<1/2$.
Estimate \eqref{eq:A-bound} then excludes a first zero of $A$ and gives
\begin{equation}\label{eq:A-positive}
A(t)>\frac12\quad(0\le t\le h),\qquad B(t)>0\quad(0<t\le h).
\end{equation}

Set
\begin{equation}\label{eq:C}
C(t)=\ell(u(t))A(t),\qquad C_h=C(h)>0.
\end{equation}
By \eqref{eq:left-variation} and \eqref{eq:ell-positive-derivative},
\[
C'=-\ell'(u)\psi(v)A-\ell(u)\psi'(v)B<0
\qquad(0<t\le h).
\]
Since $B'=v'C$ on the first half-interval,
\begin{equation}\label{eq:strict-left-inequality}
B(h)=\int_0^h C(t)v'(t)\dd t>C_hv(h).
\end{equation}

Suppose, by contradiction, that $\mu_{L,h}'(s)\le0$, and put
$m=\mu_{L,h}'(s)/\mu_{L,h}(s)\le0$.
On the second half-interval the total derivatives in \eqref{eq:AB}
satisfy
\begin{equation}\label{eq:right-variation}
A'=-\psi'(v)B,\qquad B'=v'(m+C),\qquad v'<0.
\end{equation}
On any interval starting at $h$ on which $B>0$, the function $A$
is decreasing. If $A(t)\ge0$, then $u(t)\le u(h)$ and
\eqref{eq:ell-positive-derivative} give
$C(t)=\ell(u(t))A(t)\le\ell(u(h))A(h)=C_h$.
If $A(t)<0$, the same inequality follows from $C(t)<0<C_h$.
Thus $m+C(t)\le C_h$ as long as $B$ remains positive.
Because $v'<0$, \eqref{eq:right-variation} gives $B'(t)\ge C_hv'(t)$.
Integrating yields
\begin{equation}\label{eq:no-zero-B}
B(t)\ge B(h)+C_h\bigl(v(t)-v(h)\bigr)
      =\bigl[B(h)-C_hv(h)\bigr]+C_hv(t)>0.
\end{equation}
Here \eqref{eq:strict-left-inequality} and $v(t)\ge0$ are used.
By continuity, this estimate excludes a first zero of $B$ on
$[h,2h]$. At $t=2h$, it gives $B(2h)>0$.
But the identity $v(2h;s,\mu_{L,h}(s))=0$ implies $B(2h)=0$.
The contradiction proves \eqref{eq:tail-positive}.
\end{proof}

\subsubsection{Persistence of the two critical points}

\begin{lemma}\label{lem:persistence}
Let $L_0$ be as in Proposition~\ref{prop:two-critical}, and let
$J\Subset(L_0,\infty)$ be a compact interval. For all sufficiently
small $h>0$, uniformly for $L\in J$, the function $\mu_{L,h}$
has exactly two critical points on $(0,\infty)$, denoted by
\[
0<\alpha_{L,h}<\beta_{L,h}<2.
\]
They are nondegenerate, and
\begin{equation}\label{eq:mu-sign}
\mu_{L,h}'(s)\begin{cases}
>0,&0<s<\alpha_{L,h},\\
<0,&\alpha_{L,h}<s<\beta_{L,h},\\
>0,&s>\beta_{L,h}.
\end{cases}
\end{equation}
Furthermore,
\begin{equation}\label{eq:critical-shifts}
\alpha_{L,h}=\alpha_L+O(h),\qquad
\beta_{L,h}=\beta_L+O(h),
\end{equation}
uniformly for $L\in J$.
\end{lemma}

\begin{proof}
Define
\[
\mcH(s,L,h)=
\begin{cases}
\bigl(\mu_{L,h}(s)-1\bigr)/h,&h\ne0,\\
M_L(s),&h=0,
\end{cases}
\]
using the local extension from Proposition~\ref{prop:expansion}
when $h<0$. This function is smooth on a neighborhood of
$[0,2]\times J\times\{0\}$, and
\begin{equation}\label{eq:H-C2}
\sup_{L\in J}\norm{\mcH(\cdot,L,h)-M_L}_{C^2([0,2])}=O(h).
\end{equation}

The functions $\alpha_L,\beta_L$ are continuous on $J$,
are separated from each other and from the endpoints $0,2$,
and satisfy the strict second-derivative inequalities in
\eqref{eq:M-sign}. Compactness therefore gives $\eta,c>0$
such that the intervals
\[
I_\alpha(L)=[\alpha_L-\eta,\alpha_L+\eta],\qquad
I_\beta(L)=[\beta_L-\eta,\beta_L+\eta]
\]
are disjoint and contained in $(0,2)$, and
\[
M_L''\le-c\text{ on }I_\alpha(L),\qquad
M_L''\ge c\text{ on }I_\beta(L),
\]
for every $L\in J$.
On the compact complement of the interiors of these two moving
intervals, $|M_L'|$ has a positive lower bound.
By \eqref{eq:H-C2}, the derivative $\mcH_s$ has the same sign as
$M_L'$ on this complement, is strictly decreasing on $I_\alpha(L)$,
and is strictly increasing on $I_\beta(L)$, once $h$ is small.
The endpoint signs on each interval are opposite. Thus $\mcH_s$
has exactly one zero in each interval and no other zeros in $[0,2]$.
Both zeros are simple.

Proposition~\ref{prop:tail} excludes every further critical point for
$s\ge2$, after decreasing the upper bound for $h$ if necessary.
This proves \eqref{eq:mu-sign} and nondegeneracy.
The implicit function theorem applied to $\mcH_s=0$ gives smooth
dependence of the critical points on $(L,h)$ up to $h=0$.
Their $h$-derivatives are uniformly bounded on a smaller compact
parameter neighborhood, which proves \eqref{eq:critical-shifts}.
\end{proof}

\subsection{Exact multiplicity and nondegeneracy}\label{sec:completion}

\begin{proof}[Proof of Theorem~\ref{thm:main}]
Choose $L_*>L_0$ and apply Lemma~\ref{lem:persistence} with
$J=[L_*,2L_*]$. Take $h_*>0$ small enough for that lemma,
Propositions~\ref{prop:expansion} and \ref{prop:tail}, and
condition \eqref{eq:left-condition} to hold uniformly on $J$.
For $(\lambda,h)\in\mcD$, set $L=\lambda h$ and define
\begin{equation}\label{eq:mu-critical}
\mu_+(\lambda,h)=\mu_{L,h}(\alpha_{L,h}),\qquad
\mu_-(\lambda,h)=\mu_{L,h}(\beta_{L,h}).
\end{equation}
The first point is a strict local maximum and the second a strict
local minimum. By \eqref{eq:balance} and \eqref{eq:mu-sign},
\[
1<\mu_-(\lambda,h)<\mu_+(\lambda,h).
\]
Proposition~\ref{prop:expansion} and Lemma~\ref{lem:infinity} give
\begin{equation}\label{eq:two-ends}
\lim_{s\downarrow0}\mu_{L,h}(s)=1,\qquad
\lim_{s\to\infty}\mu_{L,h}(s)=\infty.
\end{equation}
The three strictly monotone pieces of the shooting curve consequently
have ranges
\[
(1,\mu_+),\qquad (\mu_-,\mu_+),\qquad (\mu_-,\infty),
\]
when their turning-point endpoints are omitted.
Counting intersections with a horizontal line, and counting a turning
point only once, yields \eqref{eq:count}. Proposition~\ref{prop:global-curve}
ensures that these intersections account for all positive solutions.
The same description proves uniqueness of the pair $(\mu_-,\mu_+)$.
Smoothness of the critical values follows from the smooth critical-point
functions and the change of variables $L=\lambda h$.

It remains to relate the turning points to the linearized boundary
value problem. For fixed $(\lambda,h)$, use the shooting function
$F$ from \eqref{eq:shooting-F}. On the positive solution curve,
\begin{equation}\label{eq:F-derivative}
F_s+F_\mu\mu_{L,h}'=0,\qquad F_\mu<0.
\end{equation}
The linearized equation at a positive solution $u$ is
\begin{equation}\label{eq:linearized}
\bigl(\phi'(u')w'\bigr)'
 +\lambda a_{\mu,h}(t)g'(u)w=0,
\qquad w'(0)=w'(2h)=0.
\end{equation}
It is understood with continuous $w$ and $\phi'(u')w'$ across $h$.
Since $\phi'(u')$ is continuous and strictly positive, this is
equivalent to continuity of $w,w'$.
Every solution of the linearized equation satisfying the left
Neumann condition is a scalar multiple of the fixed-$\mu$ initial-value
variation $w=\partial_su$ with $w(0)=1$.
For this variation,
\[
\partial_sv=-\phi'(u')w',\qquad
F_s=\partial_sv(2h)=-w'(2h),
\]
since $u'(2h)=0$ and $\phi'(0)=1$.
Thus the linearized Neumann problem has a nonzero solution if and
only if $F_s=0$, or equivalently $\mu_{L,h}'=0$.
At either critical point its kernel is one-dimensional; elsewhere
it is trivial.

Let $s_*$ be one of the two critical initial values and let
$\mu_* =\mu_{L,h}(s_*)$. Differentiating the shooting identity twice
at $s_*$ gives
\begin{equation}\label{eq:fold-derivatives}
F_{ss}(s_*,\mu_*)=-F_\mu(s_*,\mu_*)\mu_{L,h}''(s_*)\ne0.
\end{equation}
Equivalently, the local solution curve has the form
\[
\mu-\mu_*=(s-s_*)^2J_*(s),\qquad
J_*(s)=\int_0^1(1-\theta)
\mu_{L,h}''\bigl(s_*+\theta(s-s_*)\bigr)\dd\theta,
\]
where $J_*(s_*)=\mu_{L,h}''(s_*)/2\ne0$.
The smooth local coordinate
$\xi=(s-s_*)\sqrt{|J_*(s)|}$ therefore puts the projection into
the normal form
\begin{equation}\label{eq:fold-normal-form}
\mu-\mu_*=\sgn\bigl(\mu_{L,h}''(s_*)\bigr)\xi^2.
\end{equation}
This is a nondegenerate fold. The parameterization by $s$ is regular
as a curve of functions because its derivative has value $1$ at
$t=0$. At each critical value, the other solution lies on a strictly
monotone part of the curve and is nondegenerate.
\end{proof}

\section{Critical-value asymptotics and numerical illustration}\label{sec:asymptotics}

The reduced function also gives the first-order locations of the two
critical parameter values.

\begin{corollary}\label{cor:critical-values}
Let $J\Subset(L_0,\infty)$ be a compact interval. For $L\in J$ and
sufficiently small $h>0$, define the critical values by
\eqref{eq:mu-critical}, using Lemma~\ref{lem:persistence}.
Uniformly for $L\in J$, as $h\downarrow0$,
\begin{equation}\label{eq:critical-value-asymptotics}
\begin{aligned}
\mu_+(L/h,h)&=1+hM_L(\alpha_L)+O(h^2),\\
\mu_-(L/h,h)&=1+hM_L(\beta_L)+O(h^2).
\end{aligned}
\end{equation}
In particular,
\begin{equation}\label{eq:width}
\mu_+(L/h,h)-\mu_-(L/h,h)
=h\bigl[M_L(\alpha_L)-M_L(\beta_L)\bigr]+O(h^2).
\end{equation}
\end{corollary}

\begin{proof}
Substitute \eqref{eq:critical-shifts} into
\eqref{eq:C2-expansion}. The functions $M_L$ and their first
derivatives are uniformly bounded on $[0,2]\times J$.
This gives \eqref{eq:critical-value-asymptotics}, and subtraction
gives \eqref{eq:width}.
\end{proof}

\begin{proposition}\label{prop:large-L-asymptotics}
Let $\kappa$ and $x_0$ be defined by $Q(\kappa)=1/2$ and
\eqref{eq:x0}, respectively. The reduced critical points satisfy
\begin{equation}\label{eq:reduced-critical-asymptotics}
\sqrt L\,\alpha_L\longrightarrow\sqrt\kappa,
\qquad
\beta_L\longrightarrow\sqrt{x_0}
\qquad(L\to\infty).
\end{equation}
Their critical values satisfy
\begin{equation}\label{eq:reduced-value-asymptotics}
\frac{M_L(\alpha_L)}{\sqrt L}
\longrightarrow\frac{2W(\kappa)}{\sqrt\kappa},
\qquad
M_L(\beta_L)\longrightarrow
\frac{1+2x_0}{\sqrt{x_0}}.
\end{equation}
\end{proposition}

\begin{proof}
In the first part of the proof of Proposition~\ref{prop:two-critical},
the zero
$k_L=L(e^{\alpha_L^2}-1)$ converges to $\kappa$.
Thus $L\alpha_L^2=L\log(1+k_L/L)\to\kappa$.
For the second zero, $\beta_L^2$ belongs to a fixed neighborhood
of $x_0$ and satisfies $R(\beta_L^2)=Q(L(e^{\beta_L^2}-1))\to0$.
Since $x_0$ is the only zero of $R$ in that neighborhood,
$\beta_L\to\sqrt{x_0}$.
Now $\ell(s)\sim2/s$ as $s\downarrow0$, so
\[
\frac{M_L(\alpha_L)}{\sqrt L}
=\frac{\ell(\alpha_L)}{\sqrt L}W(k_L)
\longrightarrow\frac{2W(\kappa)}{\sqrt\kappa}.
\]
For the second value, $W(Lg(\beta_L))\to1$, and
\[
\ell(\sqrt{x_0})
=\frac{2\sqrt{x_0}}{1-e^{-x_0}}
=\frac{1+2x_0}{\sqrt{x_0}}.
\]
This proves \eqref{eq:reduced-value-asymptotics}.
\end{proof}

\begin{remark}
The two limiting procedures above are distinct.
Corollary~\ref{cor:critical-values} is uniform for $L$ in a fixed
compact interval, whereas Proposition~\ref{prop:large-L-asymptotics}
concerns the reduced function as $L\to\infty$.
A joint limit in which $L$ is unbounded and $h\to0$ would require
additional uniform estimates.
\end{remark}

\begin{figure}[!htbp]
\centering
\includegraphics[width=0.72\textwidth]{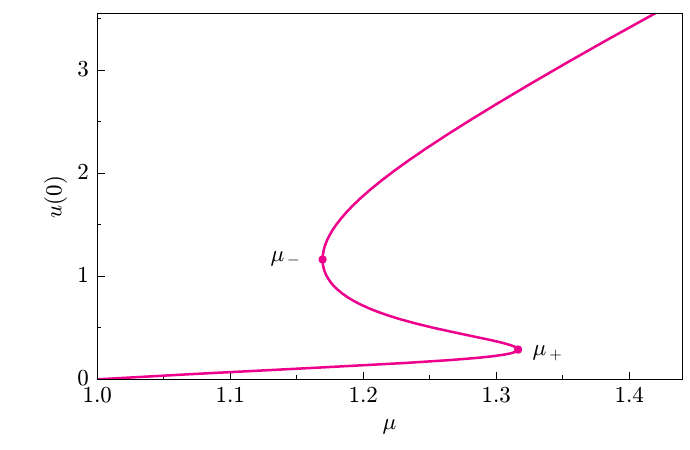}
\caption{Numerical shooting curve for \eqref{eq:problem} with
$(\lambda,h)=(400,0.05)$. The markers $\mu_-$ and $\mu_+$ indicate
the local minimum and maximum, respectively, of $\mu$ as a function
of $u(0)$.}
\label{fig:numerical-bifurcation}
\end{figure}

For numerical illustration, Figure~\ref{fig:numerical-bifurcation}
shows a portion of the shooting curve for $(\lambda,h)=(400,0.05)$.
For each sampled $s=u(0)>0$, we solve $F(s,\mu)=0$ numerically for
$\mu$ using the shooting function \eqref{eq:shooting-F}.
OpenAI Codex assisted in writing the code for this computation and the
plot. The curve was computed with SciPy and plotted with Matplotlib.
The two turning points of the computed curve lie at
$\mu_-\approx1.16942$ and $\mu_+\approx1.31644$.

\FloatBarrier

\section*{Declaration of generative AI use}

The author used OpenAI Codex to assist with the organization and
English-language editing of the manuscript, the review of references
and mathematical derivations, and the preparation of the numerical
shooting and plotting code for Figure~\ref{fig:numerical-bifurcation}.
The author reviewed and revised the material prepared with this
assistance and takes full responsibility for the content of the article.

\end{document}